\documentclass{amsart}
\usepackage{amsmath}
\usepackage{amssymb}
  \usepackage{paralist}
  \usepackage{graphics} 
  \usepackage{epsfig} 
 \usepackage[colorlinks=true]{hyperref}
\hypersetup{urlcolor=blue, citecolor=red}

\newcommand{\markattention}[1]{
#1}
\usepackage{booktabs}

\newcommand{\Branches}{\ensuremath{B}}

\newcommand{\Lottypes}{\ensuremath{L}}
\newcommand{\Sizes}{\ensuremath{S}}
\newcommand{\Pricetrajectories}{\ensuremath{T}}
\newcommand{\handlingcost}{\ensuremath{c}}

\newcommand{\addlottypecost}{\ensuremath{\delta}}
\newcommand{\testaddlottypecost}{\ensuremath{\tilde{\delta}}}
\newcommand{\markdowncost}{\ensuremath{\mu}}
\newcommand{\testmarkdowncost}{\ensuremath{\tilde{\mu}}}
\newcommand{\discount}{\ensuremath{\rho}}
\newcommand{\demand}{\ensuremath{d}}
\newcommand{\maxlots}{\ensuremath{\kappa}}
\newcommand{\kmax}{\ensuremath{k_{\max}}}
\newcommand{\kobserv}{\ensuremath{k_{\text{obs}}}}

\newcommand{\pmax}{\ensuremath{p_{\max}}}

\newcommand{\prob}{\ensuremath{\operatorname{Prob}}}

\newcommand{\lotselect}[1]{\ensuremath{x_{#1}}}

\newcommand{\lottaken}[1]{\ensuremath{y_{#1}}}
\newcommand{\lotcount}[1]{\ensuremath{z_{#1}}}

\newcommand{\stock}[2]{\ensuremath{v^{#2}_{#1}}}
\newcommand{\sales}[2]{\ensuremath{w^{#2}_{#1}}}
\newcommand{\yield}[2]{\ensuremath{r^{#2}_{#1}}}
\newcommand{\realizedyield}[1]{\ensuremath{\hat{r}_{#1}}}
\newcommand{\price}[2]{\ensuremath{u^{#2}_{#1}}}

\newcommand{\markdownused}[2]{\ensuremath{n^{#2}_{#1}}}
\newcommand{\realizedmarkdownused}[1]{\ensuremath{\hat{n}_{#1}}}
\newcommand{\MarkdownScenario}[1]{\ensuremath{a_{#1}}}

\newcommand{\Periods}{\ensuremath{K}}
\newcommand{\Prices}{\ensuremath{P}}
\newcommand{\Multiplicities}{\ensuremath{M}}
\newcommand{\Scenarios}{\ensuremath{E}}

\newcommand{\initialStock}{\ensuremath{I_{b,s}}}

\newcommand{\rro}{\ensuremath{RRO}}
\newcommand{\rrotest}{\rro_{\text{test}}}
\newcommand{\rrocontrol}{\rro_{\text{control}}}

\theoremstyle{definition}

\newtheorem{remark}{Remark}

\title[The Integrated Size and Price Optimization Problem]
      {The Integrated Size and Price Optimization Problem}

\author[M.~Kie{\ss}ling, S.~Kurz, and J. Rambau]{}

 \email{miriam.kiessling@uni-bayreuth.de}
 \email{sascha.kurz@uni-bayreuth.de}
 \email{joerg.rambau@uni-bayreuth.de}

\thanks{Supported by a grant of the \emph{Bayerische
  Forschungsstiftung, Az-727-06}}

\begin{document}
\maketitle

\centerline{\scshape Miriam Kie{\ss}ling}
\medskip
{\footnotesize
} 

\medskip

\centerline{\scshape Sascha Kurz and J\"org Rambau}
\medskip
{\footnotesize
 \centerline{Universit\"at Bayreuth,  95440 Bayreuth, Germany}
}

\bigskip


\begin{abstract}
  We present the Integrated Size and Price Optimization Problem (ISPO)
  for a fashion discounter with many branches.  Based on a two-stage
  stochastic programming model with recourse, we develop an exact
  algorithm and a production-compliant heuristic that produces small
  optimality gaps.  In a field study we show that a distribution of
  supply over branches and sizes based on ISPO solutions is
  significantly better than a one-stage optimization of the
  distribution ignoring the possibility of optimal pricing.
  
  \medskip
  
\end{abstract}

\section{Introduction}
\label{sec_introduction}

We want to decide on the one-time supply of seasonal apparel for a
real-world fashion discounter with many branches in such a way that
the expected revenue during a sales process with possible mark-downs
is maximal.  The supply process and the sales process of our partner
have several special features.  Most important among them is the
following.  The supply for a branch in a size can not be decided
independently from the other supplies: ordering and delivery are based
on pre-packed size-assortments of a product, the so-called
\emph{lots}.  A \emph{lot-type} specifies for each size the number of
items in that size in the pre-pack.  Another important restriction is
that prices can only be marked-down for a product in all branches and
all sizes, so that it is impossible to use a dynamic pricing strategy
for a product in all branches and in all sizes independently.

Moreover, most fashion products are only sold once and are never
offered again.  Thus, historical sales data can only be used on a higher
aggregation level, e.g., the average historical demand at a price in a
sales week on the \markattention{commodity group level or main commodity group} level.  Since the
average supplies per branch and size of a single product are zero,
one, or two in most cases, we can expect that historical sales data will
only give us very coarse information.

Thus, it seems reasonable to use an approach that takes into account
forecasting inaccuracies and deviations from the normal behavior.
Still we need to design a model whose stochastic parts are simple
enough to produce solutions fast while encompassing the indispensable
operational side constraints. The mere volume of merchandise handled
in a large fashion discounter (around 1000 products with, in total,
around 10 million pieces per months) requires that supply calculations
take no more than 15 minutes per product on average to keep up with
operations.  This large throughput actually implies a real-time
requirement to any algorithm used in such an environment.

In this paper we try to find a balance between modeling accuracy and
real-time compliance by incorporating only a few overall-success
scenarios of a product into the model, but stick to point forecasts
with respect to all other variabilities, like varying demands in the
branches or sizes compared to each other.

\subsection{Related Work}
\label{sec:related-work}

Linking of inventory and dynamic pricing decisions has been attacked
in \cite{Bertsimas,Swann,0979.90004,1116.90009}. More recent
approaches consider robustness considerations
\cite{springerlink:10.1007/s10107-005-0681-5} or game theoretic
aspects, like competition and equilibria \cite{adida:dynamic}.
Common to those results is the optimal control approach via fluid
approximation and/or the dynamic programming approach.  The real-world
settings of companies usually involves additional side-constraints (in
our case: the restriction on the number of used lot-types) and costs
(in our case: lot-type handling and opening costs) that would lead to
the violation of important assumptions in optimal control and that
would require very large state spaces in dynamic programming.

Dynamic pricing is a well-studied problem in the revenue management
literature (see, e.g., \cite{BitranSummer_2003,
  Gallego:1994:ODP:195519.195528,0889.90052, 1126.90372,20001212} as
examples).  Again, complicated operational side-constraints are
usually neglected in favor of a more principle study of isolated
aspects.  Again, some work has been done from a game theoretic point
of view, like strategic customers (see, e.g.,
\cite{Yin:2009:OMP:1594660.1594672}).

Classical inventory management research is less related to our topic,
since there most policies deal with the optimal way to replenish 
stock.  In our environment, no replenishment is possible.

Our first steps in capturing the operational side constraints posed by
the lot-based supply
\cite{p_median,Kiessling+Kurz+Rambau:LDP-CG-Preprint:2012} did not
take pricing into account, but estimated the consequences of inventory
decisions by a distance measure between supply and an estimated
demand.  The resulting stochastic lot-type design problem (without
reference to pricing) will serve as a template for our model of the
size optimization stage.  Since the number of possible lot-types can
be very large, we developed a branch-and-price algorithm in
\cite{Kiessling+Kurz+Rambau:LDP-CG-Preprint:2012}.  In this paper,
however, we restrict ourselves to a managable number of applicable
lot-types.

We have found no published research on field studies that apply any of
the theoretical results on dynamic pricing and/or inventory control to
a real-world environment and analyze the influence of the method in a
planned, controlled experiment.

\subsection{Our contribution}
\label{sec:our-contribution}

We present an inventory and dynamic pricing problem of a real-world
fashion discounter with a set of operational side constraints that has
been unstudied so far.  For this problem, we contribute
\begin{itemize}
\item a new model;
\item an exact branch-and-bound algorithm for benchmarking;
\item a fast heuristic for production use;
\item computational results from a field experiment with a robust
  assessment of statistical significance.
\end{itemize}
Our combined inventory and dynamic pricing problem is new compared to
already studied problems because of a combination of the following
aspects:
\begin{itemize}
\item Seasonal items.
\item No replenishing.
\item No left-over items (everything is sold at some price).
\item Early stock-outs are possible in some branches while others
  still have stock.
\item Supply must be ordered and delivered in terms of lots of a
  limited number of lot-types (pre-packed assortments of sizes of a
  single product); thus, the possible inventory decisions in each
  branch and size are severely restricted.
\item Prices must be marked-down consistently in all branches and
  sizes.
\item Prices must be taken from a small set of possible prices.
\item There are costs for handling lot-types, marginal costs for using
  another lot-type, and mark-down costs.
\end{itemize}
For the first time, we model this \emph{Integrated Size and Price
  Optimization Problem (ISPO)} as a two-stage stochastic programming
problem with recourse.  We consider mark-down decisions as recourse
actions for small success of the product.  And the real profit of a
distribution of goods over branches and sizes depends on the success
of the product: a branch and a size that receives too few items
compared to other branches and sizes produces high opportunity costs
in a high success scenario but not in a small success scenario -- too
few items on average can then become exactly right because of low
demand.  Conversely, a branch and a size that receives too many items
compared to other branches produces high mark-down losses in a small
success scenario but not in a high success scenario.

The two-stage setting is an approximation of a multi-stage setting:
ISPO is a model to decide on inventories. It ignores the flow of
information about success during the sales process (open-loop).  This
is acceptable because the price optimization stage in ISPO is only
meant to be an estimation of the recourse cost induced by the
inventory decisions.  When the sales process is actually running, we
plan to use the price optimization stage of ISPO with a rolling
horizon to utilize weekly sales information.  Thus, ISPO is
\begin{itemize}
\item less accurate than existing models in the price optimization
  stage because it estimates expected profit by an open-loop price assignment
\item much more accurate than existing models in the size optimization
  stage because it takes into account many operational side constraints.
\end{itemize}
As a contribution to modeling real-world problems, we introduce the
extended mixed integer linear programming formulation for a two-stage
stochastic programming model of ISPO.  For this model we design two
algorithms: one exact branch-and-bound method and one heuristic
method: the ping-pong heuristic.  Although both methods are kind of
tailor-made for our real-world problem, the underlying ideas could be
used in other contexts.  We state some characteristics of problems
that could be attacked by similar ideas in
Remarks~\ref{rem:relevant-structure-1}
and~\ref{rem:relevant-structure-2}.

Moreover, we performed a real-world field-study as a controlled
statistical experiment (similar to a clinical study).  We used in
parallel an existing optimization method (``old'' method) on a set of
control branches and our size optimization method based on the ISPO
model (``new'' method) on a set of test branches.  Whether a branch
was assigned to be a test or a control branch was decided randomly.
From this study we derived that in a five-month period we could
increase the mean relative realized revenue (the mean of the total
revenue divided by the maximally possible revenue) by around two
percentage points (resp.{} more than one percentage point when only a
small set of heavily cleaned up data is considered).  This means big
money in economies of scale.

The advantage of the controlled test set-up yields more: By using
robust ranking statistics exploiting the design of the experiment, we
can state that it is very unlikely (around 4\,\% probability) that
these improvements happened by chance -- and this with no assumptions
on the error distributions.  We have not seen any published results
that investigate the significance of practical results by this (or any
other) statistical method, and we consider the introduction of
controlled statistical experiments into the field of retail revenue
management as a contribution in its own right.

\emph{Nota bene:} The ``old'' method with which we compared our
``new'' method (from this paper) is \emph{not} the historical manual
solution developed at our partner's but already a one-stage size
optimization method based on the concepts in \cite{p_median}; this
``old'' method was adopted by our partner immediately after we
developed it, since it yielded obvious benefits compared to the
previous, manual solution.  The results from our field study,
therefore, evaluate the benefit of using a stochastic two-stage model
with a purely monetary objective function as opposed to a
deterministic one-stage model with a non-monetary objective function
based on a distance between supply and forecasted demand.

\subsection{Outline of the paper}
\label{sec:outline}

We state the ISPO formally in
Section~\ref{sec_formal_problem_statement}. Section~\ref{sec:modelling}
describes the extended MILP formulation of a two-stage stochastic
program with recourse that we use to model the ISPO.  In
Section~\ref{sec_algorithms} we present two algorithms: one exact
branch-and-bound solver of the MILP and the fast ping-pong
heuristic. In Section~\ref{sec_field_study} we outline the setup
of our real-world field-study. 
Section~\ref{sec_results} is devoted to computational
results: one part underpins that ping-pong can solve real-world
instances of ISPO fast with tiny optimality gaps and the other part
shows the impact of using ISPO in practice.  We conclude in
Section~\ref{sec_conclusion}.  

\section{Formal problem statement}
\label{sec_formal_problem_statement}

\noindent
We consider the distribution of supply over branches and sizes for a
single fashion article as a two-stage optimization problem.  In the
\emph{size optimization} stage (Section~\ref{sec:size-optim-stage}),
we essentially decide on a \emph{lot-type design} (see
\cite{p_median,Kiessling+Kurz+Rambau:LDP-CG-Preprint:2012}).  In the
\emph{price optimization} stage (Section~\ref{sec:price-optim-stage})
we decide on price mark-downs during the sales process depending on
the inventory induced by the first stage decisions and the overall
success of the article observed after the first sales period.

We want to maximize expected profit, where profit is given by the
yield during the price optimization stage minus costs for various
actions in both stages (Section~\ref{sec:two-stage-objective}).

\subsection{The size optimization stage}
\label{sec:size-optim-stage}

\textit{Data:} Let $\mathcal{B}$ be the set of branches of our fashion
discounter.  Let $\mathcal{L}$ be a set of applicable
\emph{lot-types}: For a set of utilized sizes $\mathcal{S}$, a
\emph{lot-type} is a vector
$(l_s)_{s\in\mathcal{S}}\in\mathbb{N}^{|\mathcal{S}|}$; it specifies the
number of pieces of each size in any pre-packed lot of this lot-type.

There is an upper bound $\overline{I}$ and a lower bound
$\underline{I}$ given on the total supply of the product over all
branches and sizes.  Moreover, there is an upper bound $\maxlots \in
\mathbb{N}$ on the number of lot-types used.\footnote{Typical ranges
  of the data are
  $|\mathcal{B}|\approx1000$,$|\mathcal{L}|\approx1000$,$3\leq|S|\leq
  7$, and $2\leq\maxlots\leq5$.  The overall bounds $\underline{I}$
  and $\overline{I}$ typically amount to around $10\,000$. We usually
  have $\underline{I}\approx\overline{I}-100$. The resulting small
  variability for the overall supply by-passes number-theoretic
  problems that may occur because the total supply can only be
  realized as a sum of cardinalities of selected lot-types.}


\textit{Decisions:} Consider an assignment of a unique lot-type $l(b)
\in \mathcal{L}$ and an assignment of a unique multiplicity $m(b)$ to
each branch~$b \in \mathcal{B}$.  These decisions specify that $m(b)$
lots of lot-type $l(b)$ are to be delivered to branch~$b$.

\textit{Decision-dependent entities:} The \emph{inventory of branch
  $b$ in size~$s$} given assignments $l(b)$ and $m(b)$ is given by
$I_{b,s}(l,m) = m(b) l(b)_s$. Moreover, the \emph{total supply}
resulting from $l(b)$ and $m(b)$ is given by $I(l,m) = \sum_{b \in
  \mathcal{B}}\sum_{s \in \mathcal{S}} I_{b,s}(l,m)$.

\subsection{The price optimization stage}
\label{sec:price-optim-stage}

\textit{Data:} We are given a supply $\initialStock$ for each branch
$b \in \Branches$ and size $s \in \Sizes$ induced by the first-stage
decisions $l$ and~$m$.  Let $k \in \Periods = \{ 0, 1, \dots, \kmax
\}$ be the index of a period, and let $\{\pi_p \}_{p \in \Prices}$ be
the set of possible prices.  In each success scenario $e \in
\Scenarios$ we know for each price $\pi_p$, each branch $b$, and each
size~$s$ the (fractional) mean demand $\demand^e_{k,p,b,s} \in
\mathbb{R}_{\ge 0}$ for the product in Period~$k$.  Moreover, a start
price $\pi_0$ and a salvage value $\pi_{\kmax}$ are given.

\textit{Realization of the demand process:} The realization of the
success scenarios takes place at the end of the $\kobserv$th period.
In all periods $0, 1, 2, \dots, \kobserv-1$ with yet unknown scenario
the start price has to be used. Since in all periods with a choice of
a price we know the success scenario we are in, only the inventory
decisions of the first stage are non-anticipative.  (This models the
situation where the success of an article can be assessed quite well
after few periods of sales.)  The first $\kobserv$ periods could be
merged to one period, but then discounting gets messy.

\textit{Decisions:} For a known success scenario~$e$ we decide for
each period $k \in \Periods \setminus \{0, \kmax\}$ on a price index
$p^e(k)$, i.e., we want to sell the product for price~$\pi_{p^e(k)}$
in period~$k$ in all branches and sizes.  This decision is taken at
the beginning of period~$\kobserv$, after the realization of the
success scenario.

\textit{Decision-dependent entities:} Let $I^e_{0, b, s} :=
\initialStock$ for all $e, b, s$.  Then for each period~$k$, a
selection $p(k)$ of price indices and an initial (fractional)
\emph{mean stock level}, denoted by $I^e_{k, b, s}(p)$, in that period
induces (fractional) \emph{mean sales}, denoted by
$\textrm{sales}^e_{k, b,s}(p)$, in period~$k$, in branch~$b$ and
size~$s$ in scenario~$e$, leading to a new mean stock level in the
next period~$k+1$.

\subsection{The two-stage objective}
\label{sec:two-stage-objective}

Using $m$ lots of a lot-type $l$ in branch~$b$ incurs a specific
\emph{lot handling cost} of $c_{l, b, m}$, e.g., a picking cost
proportional to~$m$: a lot with few pieces must be used in larger
quantities and, thus, the total supply requires more picks in
total.\footnote{A reduction of the number of picks in a high-wage
  country is the reason for the lot-based delivery from a low-wage
  country in the first place.}  For the $i$th used lot-type we have to
pay a marginal \emph{lot-type opening cost} of~$\addlottypecost_i$.
The (fractional) \emph{mean yield} in period~$k$ in branch~$b$ and
size~$s$ induced by a price assignment $p(k)$ is given by
$\textrm{yield}_{k, b, s}(p) = \pi_{p(k)} \textrm{sales}^e_{k,
  b,s}(p)$.  Each change of prices incurs a cost of~$\markdowncost$.

The goal is to find first-stage decisions such that for optimal second
stage decisions in each scenario~$e$ the \emph{expected profit}, which
is the expectation of \emph{total yield} minus \emph{lot handling
  costs} minus \emph{lot-type opening costs} minus \emph{mark-down
  costs}, is maximal.

We call this two-stage stochastic optimization problem with recourse
the \emph{Integrated Size and Price Optimization Problem (ISPO)}.

\begin{remark}\label{rem:fractionalstock}
  Fractional inventories, sales, and demands are interpreted as
  approximations of expected inventories.  In principle we could use
  many (integral) demand scenarios and integral inventory
  book-keeping.  However, the number of necessary scenarios for such a
  model would be enormous.  Thus, our scenarios model only the
  variability of the total demand induced by the overall success of
  the article. They do not model the variability of the actual demand
  with respect to periods, branches, and sizes compared to each other.
  These variabilities are ignored by using fractional values
  representing approximations of expected values.  In those cases we
  speak of \emph{mean values} rather than \emph{expected values} in
  order to distinguish the expected values that are represented by
  fractional numbers from expectations that we compute explicitly over
  all success scenarios.
\end{remark}

\section{Modelling}
\label{sec:modelling}

In the following, we develop an ILP formulation of the deterministic
equivalent of ISPO in extended form.

For the first stage (SOP) we use binary assignment variables
$\lotselect{b, \ell, m}$ to encode the independent assignment
decisions $l(b) = \ell$ and $m(b) = m$.  In the second stage we
introduce binary assignment variables for the independent second stage
assignment decision $p^e(k)$.  In order to account for the profit and
the cost, we need some more dependent variables.  We list the complete
model before we comment on the details.
\allowdisplaybreaks
\begin{align}
  \max - \sum_{b \in \Branches} \sum_{\ell \in
    \Lottypes} \sum_{m \in \Multiplicities} \lotselect{b,\ell,m} \cdot
  \handlingcost_{b,\ell,m} - \sum_{i=1}^{\maxlots} \addlottypecost_i \cdot \lotcount{i} \label{eq:ISPO:SOP-objective}\\
  + \sum_{e \in \Scenarios} \prob(e) \sum_{k \in
    \Periods} \exp(-\discount k) \Bigl( \sum_{b \in \Branches} & \sum_{s \in \Sizes} \yield{k,b,s}{e}
  - \markdowncost_k \markdownused{k}{e} \Bigr) \label{eq:ISPO:POP-objective}
\end{align}
Size Optimization Stage (SOP):
\begin{align}
  && \sum_{\ell \in \Lottypes} \sum_{m \in M} \lotselect{b,\ell,m} &= 1\ \ \ {\forall b} \label{eq:ISPO:SOP-assignment}\\
  && \sum_{m \in M} \lotselect{b,\ell,m} &\le \lottaken{\ell} \ \ \ {\forall b, \ell} \label{eq:ISPO:SOP-lotused}\\
  && \sum_{\ell \in L} \lottaken{\ell}  &\le \sum_{i=1}^{\maxlots} \lotcount{i} \label{eq:ISPO:SOP-lotcount} \\
  && \lotcount{i} &\le \lotcount{i-1}\ \ \  {\forall i} \label{eq:ISPO:SOP-usedlotsfirst}\\
  && I_{b,s} &= \sum_{\ell \in \Lottypes} \sum_{m \in M} m \cdot \ell_s \cdot \lotselect{b,\ell,m}\ \ \ {\forall b,s}\label{eq:ISPO:SOP-inventory}\\
  && I &= \sum_{b\in\Branches}\sum_{s \in \Sizes} I_{b,s}\label{eq:ISPO:SOP-totalinventory}\\
  && I &\in [\underline{I},\overline{I}] \label{eq:ISPO:SOP-totalinventorybounds}\\
  && \lotselect{b,\ell,m} &\in \{0,1\} \ \ \ {\forall b, \ell, m}\label{eq:ISPO:SOP-binaryassignment}\\
  && \lottaken{\ell} &\in \{0,1\} \ \ \ {\forall \ell}\label{eq:ISPO:SOP-binarylotused}\\
  && \lotcount{i} &\in \{0,1\} \ \ \ {\forall
    i}\label{eq:ISPO:SOP-binarylotcount}\\
  && I_{b,s}, I &\in \mathbb{Z} \ \ \ {\forall b, s}\label{eq:ISPO:SOP-nonnegativeinventory}.
\end{align}
Coupling via initial inventory:%
\begin{align}
  & & I_{b,s} - \stock{0,b,s}{e} &= 0\ \ \
  {\forall b, s, e}\label{eq:ISPO:SOP-POP-startinventory}
\end{align}
Price Optimization Stage (POP):
\begin{align}
  &&\sum_{p \in \Prices} \price{k,p}{e} &= 1 \ {\forall k, e }\label{eq:ISPO:POP-assignment}\\
  &&\price{k,0}{e} &= 1\ {\forall k < \kobserv, e} \label{eq:ISPO:POP-startprice}\\
  &&\price{\kmax,\pmax}{e} &= 1 \ {\forall e} \label{eq:ISPO:POP-salvageprice}\\  
  &&\price{k-1,p_1}{e} + \price{k,p_2}{e} &\leq 1 \ {\forall k, e, p_1,p_2 < p_1}\label{eq:ISPO:POP-nomarkup}\\
  &&\markdownused{k}{e} &\ge \price{k-1,p_1}{e} + \price{k,p_2}{e} - 1 \ {\forall k, e, p_1,p_2 \neq p_1}\label{eq:ISPO:POP-markdownused}\\
  &&\stock{k-1,b,s}{e} - \stock{k,b,s}{e} &= \sum_{p \in \Prices} \sales{k-1,b,s,p}{e}
  \ {\forall k, b, s,e}\label{eq:ISPO:POP-meanstockdynamics}\\
  &&\sum_{p \in \Prices} \sales{k,b,s,p}{e} &\leq \stock{k,b,s}{e} \ {\forall k, b, s, e }\label{eq:ISPO:POP-stockboundssales}\\
  &&\sales{k,b,s,p}{e} &\leq \price{k,p}{e} \cdot \demand^e_{k,p,b,s} \ {\forall k, b, s, p, e} \label{eq:ISPO:POP-demandboundssales}\\
  &&\yield{k,b,s}{e} &= \sum_{p \in \Prices} \pi_p \cdot \sales{k,b,s,p}{e} \ {\forall k, b, s, e}\label{eq:ISPO:POP-yieldfromsalesandprice} \\
  && \price{k,p}{e} &\in \{0,1\}\  {\forall k,p,e}\label{eq:ISPO:POP-binaryassignment}\\
  && \markdownused{k}{e} &\in \{0,1\}\  {\forall k,p,e}\label{eq:ISPO:POP-binarypricecounts}\\
  && \sales{k,b,s,p}{e} &\geq 0\  {\forall k,b,s,p,e}\label{eq:ISPO:POP-nonnegativesales}\\
  && \stock{k,b,s}{e} &\geq 0\  {\forall k,b,s,e}\label{eq:ISPO:POP-nonnegativestock}\\
  && \yield{k,b,s}{e} &\geq 0\  {\forall k,b,s,e}\label{eq:ISPO:POP-nonnegativeyield}
\end{align}
We first comment on the SOP stage model: We force an assignment of a
lot-type and a multiplicity to each branch by
Equation~\eqref{eq:ISPO:SOP-assignment}. In order to account for the
opening of a lot-type, we introduce lot-type variables
$\lottaken{\ell}$ indicating whether or not lot-type~$\ell$ is used at
all and lot-type count variables $\lotcount{i}$ that take value one if and
only if an $i$th new lot-type is
used. Equation~\eqref{eq:ISPO:SOP-lotused} guarantees that
$\lottaken{\ell} = 1$ whenever $\ell$ is assigned to at least one
branch~$b$.  Inequality~\eqref{eq:ISPO:SOP-lotcount} implies that no
more than $\maxlots$ lot-types are used.
Inequality~\eqref{eq:ISPO:SOP-usedlotsfirst} enforces that
$\lotcount{i} = 1$ implies that the number of used lot-types is at
least~$i$.  We use another dependent variable $I_{b,s}$ for the
inventory in branch~$b$ and size~$s$, and
Equation~\eqref{eq:ISPO:SOP-inventory} links this variable to the
assignment decisions.  The total inventory is then given by yet
another dependent variable $I$, computed by
Equation~\eqref{eq:ISPO:SOP-totalinventory} and enforced to stay
inside the given bounds by
Inequalities~\eqref{eq:ISPO:SOP-totalinventorybounds}.  All
independent variables have to be binary, see
\eqref{eq:ISPO:SOP-binaryassignment}
through~\eqref{eq:ISPO:SOP-binarylotcount}, while the dependent
inventory variables are integer
\eqref{eq:ISPO:SOP-nonnegativeinventory}.

Next, let us have a look at the POP stage model that is linked via the
start inventories $I_{b,s}$ to the SOP stage by
Equations~\eqref{eq:ISPO:SOP-POP-startinventory}.
Equation~\eqref{eq:ISPO:POP-assignment} enforces the assignment of
exactly one price to each period in each scenario.  That the start
price and the salvage value are fixed is expressed by
Equations~\eqref{eq:ISPO:POP-startprice}
and~\eqref{eq:ISPO:POP-salvageprice}. We forbid increasing prices by
Equation~\eqref{eq:ISPO:POP-nomarkup}.  A mark-down in period~$k$ is
indicated in the dependent binary variable $\markdownused{k}{e}$,
which is forced to one in Inequality~\eqref{eq:ISPO:POP-markdownused}
if the price has changed compared to the previous period. The
following restrictions model the dynamics of the sales process using
some dependent variables.  The fractional variable $\stock{k,b,s}{e}$
approximates the mean stock level in period~$k$ in branch~$b$ and
size~$s$ in scenario~$e$.  The fractional variable
$\sales{k,b,s,p}{e}$ measures the mean sales in period~$k$ in
branch~$b$ and size~$s$ for price~$p$ in scenario~$e$.  And
$\yield{k,b,s}{e}$ measures the mean yield in period~$k$ in branch~$b$
and size~$s$ in scenario~$e$.  (See Remark~\ref{rem:fractionalstock}
for the reason why we use fractional variables here.)
Equation~\eqref{eq:ISPO:POP-meanstockdynamics} describes the change of
stock levels from one period to
another. Inequality~\eqref{eq:ISPO:POP-stockboundssales} models that
there can be no more sales than stock, and in
Inequality~\eqref{eq:ISPO:POP-demandboundssales} we require that, only
if price~$p$ is chosen, there can be sales at price~$p$ of at most the
demand at price~$p$.  

Because the objective favors larger sales, the sales variables at a
price in an optimal solution will become exactly the minimum of stock
and demand at that price. On the level of mean values this
overestimates the mean sales; thus, this yields only an approximation.
Finally, we compute by
Equation~\eqref{eq:ISPO:POP-yieldfromsalesandprice} the yield in terms
of money.  In this POP stage, only the independent price assignment
variables need to be binary \eqref{eq:ISPO:POP-binaryassignment}.  The
dependent variables capturing the dynamics of mean stocks, sales, and
yields are required to be nonnegative in
\eqref{eq:ISPO:POP-nonnegativesales} through
\eqref{eq:ISPO:POP-nonnegativeyield}.

The objective function subtracts the costs for the handling of $m$
lots of type~$\ell$ in branch~$b$ and the lot-type opening costs for
using the first, second, \ldots, $i$th new lot-type
\eqref{eq:ISPO:SOP-objective} from the expected discounted mean yields
minus the expected discounted costs for mark-downs
\eqref{eq:ISPO:POP-objective}.

This ILP model for the deterministic equivalent of ISPO -- though yet
an approximation -- encompasses many real-world restrictions and cost
factors.  Therefore, it comes as no surprise that the branch-and-bound
phase of standard solvers (\texttt{cplex}\footnote{IBM ILOG CPLEX version 12.1}, \texttt{scip}\footnote{SCIP version 2.1.0}) did not
make any progress for months in all of our real-world instances.  The
usual real-world scale instance has $1500$ branches, $5$~sizes, some
$2000$ lot types out of which at most $5$ can be used, $4$~prices,
optimized over a time horizon of $13$ periods (usually weeks) with
respect to the expectation over $3$~success scenarios (success
above/around/below average).  And this generates a large and
complicated ILP that cannot be solved by commercial-of-the-shelf
methods at the time being.  Thus, in the next section, we present an
exact algorithm (quite fast, though not fast enough for daily
operation) and a heuristic (fast enough for daily production use, and
in all real-world tests with only tiny optimality gaps).

\section{Algorithms for the ISPO}
\label{sec_algorithms}

Since the ILP formulations presented in Section~\ref{sec:modelling} cannot be
solved directly we present an exact branch-and-bound algorithm in this section.
The main idea is to branch on the decisions of the price optimization stage.
Since the variables for the price optimization stage, see Section~\ref{sec:modelling},
are too finely grained, we consider more widescale decisions. A natural idea
is to condense the mark-down decisions in each time period to an entire price trajectory
for a given scenario $e\in\Scenarios$.

We can encode the feasible mark-down strategies or price trajectories by inserting
$\pmax-1$ symbols for a mark-down, like e.g.\ $\star$, into the sequence $1,\dots,\kmax-1$
(in period~$0$ price $\pi_0$ is fixed). An example is given by
$$
  \Big(1,2,\star,3,4,5,6,7,\star,\star,8,9,10,11,12,\star\Big),
$$
meaning that we reduce prices one time before sales period~$3$ and two times before
sales period~$8$. The fourth possible reduction is delayed after the end of the whole
sales period.  To be more precise, the concrete prices in the different sales periods
are given by:
\begin{center}
  \begin{tabular}{rrrrrrrrrrrrrrr}
    \toprule
    sales period & 0 & 1 & 2 & 3 & 4 & 5 & 6 & 7 & 8 & 9 & 10 & 11 &
    12 & 13 \\
    \midrule
    price & $\pi_0$ & $\pi_0$ & $\pi_0$ & $\pi_1$ & $\pi_1$ & $\pi_1$ & $\pi_1$ & $\pi_1$ & $\pi_3$ & $\pi_3$ & $\pi_3$ &
    $\pi_3$ & $\pi_3$ & $\pi_5$\\
    \bottomrule
  \end{tabular}
\end{center}
Having this encoding at hand, we can state that there are exactly
$$
  {|\Periods|+|\Prices|-4\choose |\Prices|-2}={\kmax+\pmax-2\choose \pmax-1}
$$
feasible mark-down strategies or price trajectories, i.e., in our example we have
${16 \choose 4}=1820$~price trajectories for each scenario. 

The exact details of the branching steps and the bounding step are
outlined in Subsection~\ref{sec:an-exact-branch}. As upper bounds for
the remaining size optimization stage we use both tailored
combinatorial bounds, see Subsection~\ref{sec:upper-bounds-SOP}, which
can be computed efficiently, and linear relaxations. Algorithmically
the efficient computation of those lower bounds is based on the fast
solution of a certain subproblem, see Subsection~\ref{sec:adjust}. To
remove some complications caused by these minor details one can also
assume that we solve the SOP subproblem by using the corresponding ILP
formulation directly, see Subsection~\ref{sec:exact_sol_SOP}, without
computing any cheaper bounds.

In Section~\ref{sec_results} we present computational results for the
proposed brand-and-bound algorithm.

As a heuristic, that can as well be used at the start of the
branch-and-bound, we present the so-called \emph{ping-pong} heuristic
in Subsection~\ref{sec:ping-pong-heuristic}. It will turn out that it
achieves a very good solution quality while requiring only little
computation time. The underlying idea is to iteratively solve the
separate subproblems of the size optimization and the price
optimization stage. The temporary solution of one subproblem is then
taken as an input for the other subproblem. To this end we present an
exact but rather easy, exact algorithm for the prize optimization
stage in Subsection~\ref{sec:exact_sol_POP}. To speed up ping-pong, we
can use a heuristic for the SOP from \cite{p_median} for the size
optimization subproblem, which we recall in Subsection~\ref{sec:SFA}.

The remaining part of this section is arranged as follows. At first we
present our workhorses for the solution of intrinsic subproblems in
Subsections~\ref{sec:adjust}--\ref{sec:exact_sol_POP}. For a first
reading these can be skipped.  In Subsection~\ref{sec:an-exact-branch}
we present our main branch-and-bound algorithm and in
Subsection~\ref{sec:ping-pong-heuristic} the ping-pong heuristic.

\subsection{Workhorse 1: Adjusting supplies to the total-supply
  constraints at minimal cost}
\label{sec:adjust}
Suppose that we want to solve the following rather general binary
linear problem:

\begin{align}
  && \min \sum_{v\in\mathcal{V}}\sum_{a\in \mathcal{A}} \sum_{b \in \mathcal{B}} \psi(v,a,b)\cdot x_{v,a,b}\nonumber\\
  && \sum_{a\in \mathcal{A}} \sum_{b \in \mathcal{B}} x_{v,a,b} &= 1\ \ \ {\forall v\in\mathcal{V}} \\
  && \sum_{v\in\mathcal{V}}\sum_{a\in\mathcal{A}}\sum_{b\in\mathcal{B}} \varphi(a,b)\cdot x_{v,a,b} &\in [\underline{R},\overline{R}]\label{ressource_constraint}\\
  && x_{v,a,b} &\in \{0,1\} \ \ \ {\forall v\in\mathcal{V},a\in\mathcal{A},b\in\mathcal{B}},
\end{align}
where
\begin{enumerate}
 \item[(1)] $\varphi(a,b)$ is monotonously increasing in $b$ and
 \item[(2)] $\psi(v,a,b)$ is convex in $b$.
\end{enumerate}
The continuous relaxation of this problem can be solved by a greedy
approach: In the initialization phase we determine for each
$v\in\mathcal{V}$ and each $a\in\mathcal{A}$ the optimal value
$b_{v,a}\in\mathcal{B}$ with minimal costs $\psi(v,a,b)$ using binary
search.  This can be done in
$O\bigl(|\mathcal{V}|\cdot|\mathcal{A}|\cdot\log(|\mathcal{B}|)\bigr)$
steps.  By $v(a)$ we denote that element $a\in\mathcal{A}$ that
minimizes $\psi(v,a,b_{v,a})$ and by $v(b)$ the corresponding value
$b_{v,v(a)}$. With this we set $x_{v,v(a),v(b)}=1$ for all
$v\in\mathcal{V}$.  All other values are set to zero. If
Inequality~(\ref{ressource_constraint}) is satisfied by pure chance,
then the current assignment of the $x_{v,a,b}$-variables yields a
globally optimal solution.

Otherwise we have to adjust the $x_{v,a,b}$ in order to satisfy the
resource constraint. For brevity we only discuss the case where
$\sum_{v\in\mathcal{V}}\sum_{a\in\mathcal{A}}\sum_{b\in\mathcal{B}}
\varphi(a,b)\cdot x_{v,a,b}>\overline{R}$.  The other case is
analogous.  Here, we iteratively have to take away resources from some
of the $v\in \mathcal{V}$. To this end, we introduce \emph{relative
  costs $\Delta_{v,a}^{-}$} for each $v\in\mathcal{V}$ and each
alternative $a\in\mathcal{A}$. Using
\begin{equation}
  \beta_v(a)=\max\{b\in\mathcal{B}\mid \varphi(a,b)<\varphi(v(a),v(b))\},
\end{equation}
for a given $a\in\mathcal{A}$, the relative costs for changing the
pair $\bigl(v(a),v(b)\bigr)$ to $\bigl(a,\beta_v(a)\bigr)$ are given
by
\begin{equation}
  \Delta_{v,a}^{-}=\frac{\psi\bigl(v,v(a),v(b)\bigr)-\psi\bigl(v,a,\beta_v(a)\bigr)}
  {\varphi\bigl(v(a),v(b)\bigr)-\varphi\bigl(a,\beta_v(a)\bigr)}
\end{equation}
per resource item. If $\{b\in\mathcal{B}\mid
\varphi(a,b)<\varphi(v(a),v(b))\}=\emptyset$ we set the corresponding relative
costs to $\Delta_{v,a}^{-}=\infty$. By $\omega(v)\in \mathcal{A}$ we
denote the alternative with smallest relative costs
$\Delta_v^{-1}:=\min\{\Delta_{v,a}^{-}\mid
a\in\mathcal{A}\}=\Delta_{v,\omega(v)}^{-}$ and by $v^\star$ we denote
the element in $\mathcal{V}$, where the pair
$\bigl(\omega(v),\beta_v(\omega(v))\bigr)$ attains the globally
smallest relative costs.

As abbreviations we use $R=\sum_{v\in\mathcal{V}} \varphi(v(a),v(b))$
and $\delta=\varphi\bigl(v^\star(a),v^\star(b)\bigr)
-\varphi\bigl(\omega(v^\star),\beta_v(\omega(v^\star))\bigr)$. Due to
the convexity of the target function we can state the following:
\begin{enumerate}[(a)]
\item If $\Delta_{v^\star}^-=\infty$, then the problem is
  infeasible.
\item If $R-\delta\ge \overline{R}$, then after performing the
  greedily optimal replacement of the pair
  $\bigl(v^\star(a),v^\star(b)\bigr)$ by
  $\bigl(\omega(v^\star),\beta_v(\omega(v^\star))\bigr)$ the new
  assignments correspond to an optimal solution of our problem, where
  Inequality~(\ref{ressource_constraint}) is replaced by
  $\sum_{v\in\mathcal{V}}\sum_{a\in\mathcal{A}}\sum_{b\in\mathcal{B}}
  \varphi(a,b)\cdot x_{v,a,b}\le R-\delta$.
\item If $R-\delta<\overline{R}$ we obtain the optimal solution
  of our relaxed problem with fractional variables $x_{v,a,b}$ by utilizing a
  suitable linear combination of the old and the new assignment.
\end{enumerate}
Thus, after a finite number of iterations, depending at most linearly
on the difference between the initial overall resource consumption and
$\overline{R}$, we obtain the optimal solution of the problem with at
most two fractional variables~$x_{v,a,b}$.

We remark that is also possible to solve the integral problem by
utilizing a branch-and-bound approach -- we do not go into the details here.

\subsection{Workhorse 2: Upper bounds for the Size Optimization Problem}
\label{sec:upper-bounds-SOP}

In later parts of the algorithms we need a computationally cheap dual,
i.e., upper bound for the SOP with a fixed scenario~$e$ and price
trajectory~$t$. We establish our first upper bound based on the
integrality of the individual supply for each branch in each size,
relaxing the constraints arising from a lot-based distribution.

If we supply branch~$b$ in size~$s$ with $I_{b,s}$ items, then we can
compute the costs $\lambda_{b,s}^{e,t}(I_{b,s}):=\sum_{k \in\Periods}
\exp(-\discount k)\yield{k,b,s}{e}\in\mathbb{R}_{\ge0}$ directly using
$e$, $t$, and $I_{b,s}$ to evaluate the dependent
variables~$\yield{k,b,s}{e}$.

The supply of branch~$b$ with lot-type~$l$ in multiplicity~$m$ results
in handling costs of $c_{b,l,m}$. Let $\tilde{c}_{b,s}(i)\ge 0$ be
that part of the costs that can be associated with a supply of
branch~$b$ in size~$s$ with $i$~items.

By $\hat{\lambda}_{b,s}^{e,t}$ we denote the maximum value of
$\lambda_{b,s}^{e,t}(I_{b,s})-\tilde{c}_{b,s}(I_{b,s})$ for all
achievable supplies $I_{b,s}$.  We comment only briefly on how to
compute these values fast: The most simple thing that always works, is
to exhaustively enumerate the set of possible $I_{b,s}$ (if we assume
it to be finite). Once $\yield{k,b,s}{e}$ and $-\tilde{c}_{b,s}$ are
concave functions in $I_{b,s}$ we can more sophistically compute the
maximum using nested intervals.

As we have to use at least one lot-type and we assume
$\addlottypecost_i\ge 0$, the costs of the objective function of ISPO
that can be associated with trajectory~$t$ are bounded from above by
\begin{equation}
  \sum_{e \in \Scenarios} \prob(e)\Bigl(-\addlottypecost_1+\sum_{b\in\Branches}\sum_{s\in\Sizes}
  \hat{\lambda}_{b,s}^{e,t}\Bigr).
\end{equation}
Using the general method from Subsection~\ref{sec:adjust},
we can additionally incorporate the restrictions on the overall supply
(where we assume that the convexity condition is satisfied, which is
the case in our setting). Here, $\mathcal{V}$ are the pairs $(b,s)$ of
branches and sizes, $\mathcal{A}$ is an used set consisting of one
element, and $\mathcal{B}$ is the set of possible supplies $I_{b,s}$
to a branch~$b$ in size~$s$.

Another possibility to further tighten the upper bound is to
incorporate the fact that the branches have to be supplied using
lot-types in a certain multiplicity. So, in an initialization phase
one can compute a locally best-fitting lot-type and multiplicity for
each branch separately. If the number of lot-types and possible
multiplicities is small enough, then this can be done simply by
exhaustive enumeration. For more sophisticated methods based on a
suitable parameterization of the set of applicable lot-types we refer
to~\cite{Kiessling+Kurz+Rambau:LDP-CG-Preprint:2012}. The restrictions on the
overall supply can then be incorporated by using the algorithm from
Subsection~\ref{sec:adjust}, where $\mathcal{V}$ is the set
of branches $\Branches$, $\mathcal{A}$ is the set of lot-types
$\Lottypes$, and $\mathcal{B}$ is the set of multiplicities
$\Multiplicities$. In other words, we have relaxed the restriction to
a certain number of used different lot-types and ignored the
corresponding costs.

In our concrete application we have used all three mentioned upper
bounds. Thus, our computational results rely on convexity; in all
other cases the algorithm has to use the first bound only and will
usually be slower.

\subsection{Workhorse 3: A heuristic for the Size Optimization Problem}
\label{sec:SFA}
In \cite{p_median} the so-called \emph{Score-Fix-Adjust (SFA)}
heuristic was proposed for the Lot-Type Design Problem (LDP).  The
LDP is directly related to the SOP stage of ISPO with a fixed scenario
and a fixed mark-down strategy.  In order to apply SFA to the SOP stage,
we need to modify SFA to cope with opening and the handling costs for
lot-types. Fortunately, this can be achieved by a suitable
modification of the cost coefficients in an ordinary LDP:
Incorporating the handling costs in the cost coefficients of an LDP is
simply done by adding the handling cost $\handlingcost_{b,\ell,m}$ to
the cost coefficient of~$\lotselect{b,\ell,m}$. The opening costs for
lot-types can be taken into account by solving for each possible
number of lot-types $1$ through $\maxlots$ an individual LDP with a
prescribed numbers of used lot-types, add the corresponding opening
costs to the optimal objective function, and pick the best option in
hindsight.

For completeness, we briefly describe the underlying idea of the SFA
heuristic for the LDP in the special form we use it \markattention{anytime we solve the SOP}
stage. For each branch we determine the three locally best fitting
lot-types and add a score of 100 to the best fitting lot-type, a score
of 10 to the second best fitting lot-type and a score of 1 to the
third best fitting lot-types. (Of course this can be generalized to
the first $t$ best fitting lot-types and different scoring schemes.)
With this we have implicitly assigned a score to each lot-type
$l\in\Lottypes$, where most of the lot-types obtain the score zero.
We can extend this scoring to the $k$-subsets of $\Lottypes$ by
summing up the individual scores so that we implicitly get an order of
the ${{|\Lottypes|} \choose k}$ many feasible lot-type combinations.
With this we traverse the $k$-subsets of $\Lottypes$ in descending
order, where ties are broken arbitrarily. (The crucial observation is
that this can be done without explicitly generating all such subsets
beforehand.)  In the fixing step we assume that the applicable
lot-types are restricted to the current $k$-subset of~$\Lottypes$.
Now we are in the situation where we can apply the algorithm from
Subsection~\ref{sec:adjust}.  In the initialization we start with a
locally optimal assignment of lot-types and multiplicities. We choose
$\mathcal{V}=\Branches$, $\mathcal{A}=\Lottypes$, and
$\mathcal{B}=\Multiplicities$.

We remark that the SFA-heuristic reliably produces close-to-optimal
solutions on real-world instances, see \cite{p_median}.

\subsection{Workhorse 4: Exact solution of the Size Optimization Problem}
\label{sec:exact_sol_SOP}
For a given scenario~$e$ and a given price trajectory~$t$, the ISPO is
simplified to an LDP with modified cost coefficients. This subproblem
can, e.g., be solved by utilizing the restricted version of the ILP
formulation given in Section~\ref{sec:modelling} -- and we do this for
obtaining the computational results in this paper. For more
sophisticated algorithms we refer to
\cite{Kiessling+Kurz+Rambau:LDP-CG-Preprint:2012}, where a tailored
branch-and-price algorithm is proposed that can handle millions of
lot-types.

\subsection{Workhorse 5: Exact solution of the Price Optimization Problem}
\label{sec:exact_sol_POP}
For a given scenario~$e$, a given mark-down strategy~$t$, and the given
initial supplies $I_{b,s}$ for all branches and sizes we can easily
compute the number of sold items per branch, size, and period. Since
in any reasonable setting all prices except maybe the salvage value
are positive, we conclude that in any optimal solution the number of
sold items is exactly the minimum of stock and demand in each period.
With this all other dependent variables of the ILP formulation for the
POP in Section~\ref{sec:modelling} can be computed. Therefore, we can
solve the POP stage by exhaustive enumeration of all possible mark-down
strategies. This can be done in
$O\left(\left|\Branches\right|\cdot\left|\Sizes\right|\cdot\left|\Periods\right|
  \cdot\left|\Pricetrajectories\right|\right)$ steps, which is
possible in all practical situations ($\left|\Branches\right|\approx1000$,$3\leq\left|\Sizes\right|\leq7$,$\left|\Periods\right|=13$,$\left|\Pricetrajectories\right|=1820$) we have encountered so far.

\subsection{An exact branch-and-bound algorithm}
\label{sec:an-exact-branch}
In this subsection we propose our main algorithm -- a customized
branch-and-bound algorithm.  We branch on maps ``scenario $\mapsto$
price trajectory''. A node at depth~$j$ then corresponds to all such maps
with the images of the first $j$~scenarios fixed.  The leaves are the
maps with fixed images for all scenarios.  The cost of a leaf can be
computed by solving an LDP problem from
Section~\ref{sec:exact_sol_SOP} (Workhorse~4), crucially depending on
the method of Section~\ref{sec:adjust} (Workhorse~1).  As dual bounds
we utilize the upper bounds from Subsection~\ref{sec:upper-bounds-SOP}
(Workhorse~2). As primal bounds we employ the heuristically found
solutions from Subsection~\ref{sec:SFA} (Workhorse~3), again using the
method in Secion~\ref{sec:adjust} (Workhorse~1).  In the branching
step we extend a partially defined map in a node by all possible price
trajectories for the next scenario.


In the following, we present the detailed implementation of the above
concept.  In an initialization step we compute for each scenario~$e$
and each price trajectory~$t$ a combinatorial upper bound using the
algorithms from Subsection~\ref{sec:upper-bounds-SOP}. The bound
$\Gamma(e,t)$ is saved for each pair $(e,t)$ and possibly updated
later on. Using these bounds we label the price trajectories in
ascending order: $t_1^e,\dots,t_{|\Pricetrajectories|}^e$.

Next we consider the branching step at a node of Depth~$j$, where the
price trajectories $\xi_j\in\Pricetrajectories$ of the first $j$
scenarios are already fixed. If $j<|\Scenarios|$ then we consider the
possible price trajectories for scenario $j+1$. We loop from $i=1$ to
$i=|\Pricetrajectories|$ and consider price trajectory
$t_i^{j+1}$. Now we compute the upper bound
\begin{align}
  &\sum_{h=1}^j \prob(h)\cdot \Gamma(h,\xi_h)+\prob(j+1)\cdot\Gamma(j+1,t_i^{j+1})\nonumber\\
  {} + {} &\sum_{h=j+2}^{|\Scenarios|} \prob(h)\cdot
  \max\bigl\{\Gamma(h,t)\mid t\in\Pricetrajectories\bigr\}\label{upper_bound_bab}
\end{align}
for the ISPO where the first $j+1$ price trajectories are fixed to
$\xi_h$. If this bound is smaller than the best found integral
solution of the ISPO, then we can prune all price trajectories
$t_h^{j+1}$ for $h\ge i$. Otherwise we check how the bound
$\Gamma(j+1,t_i^{j+1})$ was computed. If it was computed using the
combinatorial relaxations from Subsection~\ref{sec:upper-bounds-SOP},
then we compute the LP bound from the restricted ILP model, see
Subsection~\ref{sec:exact_sol_SOP}\footnote{Here we apply warm-start
  techniques and initialize the LP with a basis solution of a
  \textit{similar} price trajectory within the same scenario -- if
  available.}, and possibly update the bound $\Gamma(j+1,t_i^{j+1})$.
If the updated upper bound~(\ref{upper_bound_bab}) is still to weak to
prune the subtree, we fix $\xi_{j+1}=t_i^{j+1}$ and continue at the
next node.\footnote{Another possibility to improve the upper
  bound~(\ref{upper_bound_bab}) is to replace the first sum and the
  central summand by the optimal target value of the LP arising from
  ISPO restricted to the first $j+1$ scenarios -- we have not used
  this improvement in our computational results.}

In the leaves, where all price trajectories are fixed, we solve the
remaining SOP, see Subsection~\ref{sec:exact_sol_SOP}.


\subsection{The ping-pong heuristic}
\label{sec:ping-pong-heuristic}
Since the exact algorithm is still not fast enough for daily
production (see Section~\ref{sec:performance}), we have developed a
fast heuristic. The main idea is to alternatingly fix the independent
variables of one stage and compute the optimal remaining variables;
thereafter, the resulting independent variables of the other stage are
fixed, and the remaining variables are computed optimally.  And so
on. To be more precise, if the independent decisions of the first
stage, i.e., the supply of the branches with lot-types in a certain
multiplicity -- in other words: the $x_{b,l,m}$ -- are given, then one
can easily solve the prize optimization problem of the second stage by
exhaustively enumerating all possible mark-down strategies in all
scenarios separately. If for the other direction the independent
decisions of the second stage, i.e., the $\MarkdownScenario{e,t}$, are
fixed, then the remaining problem reduces to the SOP stage, which is
essentially an LDP with a modified cost function. 

The idea now is to use the (close-to-) optimal solution of one of
these two subproblems as input for the other subproblem and to iterate
this until the algorithm stays at a solution.  We hope that the
solution that is not changed anymore is a good solution.

More specifically, we perform the following steps:
\begin{enumerate}
\item \label{itm:ping-pong:1} \textit{Initialization:} In all
  scenarios we choose the mark-down strategy which produces the best
  combinatorial bound, see Section~\ref{sec:upper-bounds-SOP}
  (Workhorse~2).
\item \label{itm:ping-pong:2} Given the mark-down strategies for all
  scenarios we heuristically solve the remaining SOP stage with the
  SFA heuristic, see Section~\ref{sec:SFA} (Workhorse~3).
\item \label{itm:ping-pong:3} Given the initial supply of the branches
  we exactly solve the prize optimization problem of the second stage,
  see Section~\ref{sec:exact_sol_POP} (Workhorse~5).
\item \label{itm:ping-pong:4} As long as the solutions of
  Steps~\ref{itm:ping-pong:2} and~\ref{itm:ping-pong:3} have not
  converged and the number of iterations is below a certain threshold,
  we proceed with Step~\ref{itm:ping-pong:2}.
\end{enumerate}
Finally we output the best solution of the ISPO found in
Step~\ref{itm:ping-pong:2}.

\begin{remark}\label{rem:relevant-structure-1}
  The details of our branch-and-bound method are involved.  However,
  there is one crucial property of the problem because of which the
  method works: Our problem has a \emph{reversible two-stage
    structure}.  This means: the \emph{independent second stage}
  variables (in our case the maps from scenarios to price assignments)
  can be interpreted as \emph{independent first stage} decisions.  The
  independent first stage variables and all dependent variables can
  then be seen as second stage variables.  In our setting fixing the
  independent decision variables of one stage does not even imply any
  restrictions to the feasible set of the independent variables in the
  other stage.  We call this \emph{reversible complete recourse}. In
  general, a heuristic like ping-pong is promising if fixing the
  independent variables from one stage leaves over a rich feasible set
  for the other stage, hopefully always containing improving
  solutions.  In our case, fixing a price trajectory to a scenario
  does not influence the feasibility of supply.  This is the case for
  all inventory problems where the price dependent demand can be determined a priori.
\end{remark}

\begin{remark}\label{rem:relevant-structure-2}
  The principle of the ping-pong heuristic is similar to the principle
  of evolutionary algorithms, see for example~\cite{MIETTINEN1999}.
  The idea of evolutionary algorithms is to assign a so-called
  \emph{fitness}-function to the solutions and iteratively in a
  \emph{selection}-step to combine the best-solutions to get solutions
  with higher \emph{fitness}. This is done until convergence.  In our
  case the \emph{fitness} of the supply in terms of lot-types is given
  by the expected revenue by the price optimization stage. By
  combining the local optimal supply with the local optimal mark-down
  strategy we possibly get a supply which results in higher revenue.
  One could also connect the principle of our ping-pong heuristics
  with the principle of bilevel programming.  A bilevel program
  consists of an upper-level and a lower-level optimization
  problem. The lower-level problems considers a variable $x$ as a
  parameter to compute the optimal value of a variable $y$ while the
  upper-level problem obtains the optimal value of $x$ by using the
  value of $y$ computed in the lower-level
  problem~\cite{ConejoCastilloMinguezGarcia-Bertrand201011}.  In our
  case -- by virtue of reversible complete recourse -- we can see the
  size optimization stage and also the price optimization stage as
  both, as upper-level and as lower-level subproblems.
\end{remark}

\section{Setup of the field study -- a controlled experiment}
\label{sec_field_study}

We performed a real-world field-study as a controlled statistical
experiment. On the one side we use the method currently applied by our
business partner (named ``old'' method\footnote{The ``old'' method
  represents a lot-type optimization method that does not take into
  account the pricing stage but estimates the gain and the loss of a
  distribution of supply by a distance measure between the supply
  induced by the lot-design and the forecasted mean demand.  This is
  not the method originally used by our business partner prior to our
  cooperation.  The old method is essentially the Score-Fix-Adjust
  heuristic presented in \cite{p_median} and refined in
  \cite{Kiessling+Kurz+Rambau:LDP-CG-Preprint:2012}.  This method
  performed so much better than the manual solutions without
  optimization that it was put into operation immediately.} in the
following) on a set of control branches. We compared this supply
strategy with the results that ISPO produced (named the ``new'' method
in the following) on a set of test branches.

The field-study ran from end of May until end of September 2011 for
81~articles from three different commodity groups -- women
overgarments fashion (wof), women overgarments classic (woc) and women
underwear (wu).  It was necessary to select a subset of articles for
the field study because the orders had already been placed in terms of
lot-types, and the adaption of the supply for the test branches to the
results of the new method was a far too expensive logistic operation to
be carried out for each article.

Since for all advertized products, in particular for those in the
field study, it is obligatory to supply each branch with at least one
piece in each size, the degree of freedom in distributing the supply
is severely restricted: small branches are very likely to receive the
one-for-all-sizes lot, leaving fewer options for the larger branches,
since the total supply is essentially fixed.

The sales process for the articles in the field study started between
May 2011 and mid of June 2011 so that all articles could be observed
for a time period of 15 to 17 weeks. Some further relevant properties
of the used test articles are stated in Table~\ref{testartikel}.

\begin{table}[ht!]
  \begin{center}\sffamily\footnotesize
    \begin{tabular}{ccc}
      \toprule
      \textbf{commodity group} & \textbf{number of articles} & \textbf{number of sizes}  \\
      \midrule
      wof & 9 & 6 \\
      woc & 9 & 3 \\
      wu  & 5 & 6 \\
      \bottomrule  
    \end{tabular}
  \end{center}
 \caption{Properties of the test articles.}
 \label{testartikel}
\end{table}

In order to obtain statistically assessable results, we grouped the
branches involved in the field study into $30$ pairs according to
economic key figures, like the size of the stores and revenue.
Whether a branch was assigned to be a test or a control branch in such
a pair was then decided randomly. In Section~\ref{sec:fieldstudy} we
will benefit from this controlled test set-up and apply robust ranking
statistics without assuming anything about the underlying error
distributions.

The test branches were supplied according to close-to-optimal
solutions of ISPO computed by the ping-pong heuristic. The ISPO for
each article in the test selection was set-up for all branches and
sizes. Since there are global constraints for the overall number of
supplied items we actually \emph{computed} the supply for \emph{all}
branches with the new method -- and so did our our project partner
with the old method. Our proposed supply was then \emph{implemented}
only for the \emph{test branches}; the supply of the \emph{control
  branches} (and all remaining branches) was \emph{implemented} as
computed by the old method by our project partner.

The demand $\demand^e_{k,p,b,s}$ was estimated based on historical
sales data of articles from the same commodity group. This is highly
non-trivial, and there are no publications claiming a ``best'' method
for this important building block.  We essentially took non-parametric
estimations of average values over commodity groups and interpolated
values with too few observation linearly, which turned out to be more
reliable than parametric estimators: the corresponding assumptions
(exponentially decreasing stock, isoelasticitic dependence of the
demand on the price, \ldots) appeared to be doubtful guesses in our
environment at best.

Table~\ref{tab:parameters} shows the parameter setting we used in ISPO
for the field study.  In order not to reveal company internals, we
printed the values with respect to artificial but consistent monetary
units. Important is that the handling cost $\handlingcost_{b,\ell,m}$
contains a term linear in the multiplicity~$m$: this way there is some
incentive for the optimization to prefer lots that produce fewer picks
in the warehouse, which is the true designation of lots in the first
place. Our lot-type opening costs $\addlottypecost_i$ were estimated
on the basis of a thourough cost accounting.  This cost accounting
also revealed that more than four lot-types can only be handled if the
area for internal stock-turnover is increased substantially.  The
discount factor~$\rho$ was derived from an estimation of the capital
binding cost.  Whenever other reasons than interest rates favor faster
stock-outs this can be increased.  The fact that we did not account
for mark-down costs $\markdowncost_k$ just reflects the fact that at
the time of the design of the experiment our partner could simply not
provide a realistic value for this.  (Meanwhile, we have an estimation
for this as well.)
\begin{table}[ht!]
  \begin{center}\footnotesize\sffamily
    \begin{tabular}{rl}
      \toprule
      \textbf{parameter} & \textbf{setting} \\
      \midrule
      $\maxlots$ & 4 \\
      $\handlingcost_{b,\ell,m}$ & acquisition price + $m \cdot  0.0545$ (pick cost) \\
      $\addlottypecost_1$ & $100$ \\
      $\addlottypecost_i$, $i>1$ &  $50$ \\
      $\Scenarios$ & $\{\text{low},\text{normal},\text{high}\}$ (period-$0$ sales $\in$
      $[0,10\,\%)$/$[10\,\%,30\,\%]$/$(30\,\%,100\,\%]$)\\
      $\demand^{\text{normal}}_{k,p,b,s}$ & from empirical
      distribution and interpolation of historical sales in commodity group\\
      $\demand^e_{k,p,b,s}$ & $\alpha \times \demand^{\text{normal}}_{k,p,b,s}$ ($\alpha$ from historical~sales in scenario~$e$ compared to scenario~$\text{normal}$)\\
      $\prob(e)$ & from empirical distribution of historical
      sales in commodity group\\
      $\kobserv$ & $2$ (i.e., realization of~$e$ and earliest
      mark-down after $2$ periods)\\
      $\kmax$ & periods (= weeks) until end of season (article dependent)\\
      $\discount$ & $0.000974868$\\
      $\pmax$ & 4 (five prices including start price and salvage value)\\
      $\markdowncost_k$ & 0\\
      \bottomrule
    \end{tabular}
  \end{center}
 \caption{Parameter setting for the field study.}
 \label{tab:parameters}
\end{table}

Since ISPO computes a supply for each branch and size under the
assumptions that optimal (open loop) prices are chosen later on, we
needed to gain control over the price optimization phase as well.

We had to decide whether we should, in the test branches,
\begin{itemize}
\item use an open loop price policy based on our POP model for the
  second stage after the revelation of the success scenario (``POP''),
\item use a closed-loop pricing policy based on our POP model applied
  with receding horizon (``RH-POP''), or
\item use whatever our partner uses for marking down prices (``manual'').
\end{itemize}
Moreover, we had to decide whether or not to allow special offers and
campaigns in the test and control branches.

The ``POP'' option would be closest to ISPO as a model, but farthest
from practice: nobody would ever ignore up-to-date sales information
for the mark-down decisions.  The ``manual'' option is closest to
practice if the mark-down process is not planned to change, but it
leaves it open whether the applied mark-down policy is anywhere close
to optimal.  We chose the second option ``RH-POP'' because of two
reasons:
\begin{itemize}
\item We tested the resulting closed-loop policy (in a different field
  study), and it performed slightly favorably compared to the manual
  policy that is currently in operation.  Thus, the difference to the
  actual operations seemed not too large and it could potentially
  replace the mark-down system currently in use. Therefore, we
  concluded that it would produce results not too far from real
  operations.
\item The POP stage of ISPO can be viewed as an estimation for the
  results obtained by the RH-POP. Therefore, we concluded that RH-POP
  would be not too far from the model assumptions either.
\end{itemize}
In order to assess better the practicability of our new method, we
decided not to forbid campaigns -- triggered by external reasons like
new competing stores -- in the test and control branches: A method the
performance of which vitally relies on laboratory conditions with all
exogeneous disturbances removed cannot be used in practice anyway.

We used RH-POP in a slightly refined way:
For each non-negative real number we defined a success scenario:
Scenario ``$1.0$'' means that the overall mean demand is as forecasted
(the mean of the historical success in the commodity group).  In
general, Scenario ``$\alpha$'' means that each mean demand is actually
$\alpha$ times the forecasted mean demand.  At the end of each period,
we updated our estimation of~$\alpha$ by comparison of our predicted
demands (based on the old $\alpha$) with the demands observed in the
just completed period.  Then, we used POP to compute a new open-loop
price policy. Whenever the optimal price trajectory suggested a
mark-down in the following two time periods, we advised our industry
partner to implement exactly this mark-down.  The method is reminiscent
of model predictive control~\cite{Gruene+Pannek:NLMPC:2011}.

\newcommand{\significancel}{\alpha} 
\newcommand{\win}{w}
\newcommand{\loss}{l}
\newcommand{\samplesize}{n}

\section{Computational results}
\label{sec_results}

In this section we report on extensive computational results about
\begin{itemize}
\item the technical performance of our algorithms in the laboratory;
\item the practical performance of their solutions in the real-world
  field study.
\end{itemize}

\subsection{Performance of the exact algorithm versus ping-pong}
\label{sec:performance}

Table~\ref{tab:runtime-results} shows the performances of the exact
branch-and-bound algorithm compared to the performance of the
ping-pong heuristic.

We ran branch-and-bound and ping-pong on many real-world instances with
\begin{itemize}
\item more than $1000$ branches
\item more than $1000$ applicable lot-types, out of which at most $5$ can
  be used in a lot-type design
\item $13$ periods $0, \ldots, 12$
\item $4$ prices that can be non-increasingly set in periods $1$
  through~$11$.
\item $3$ scenarios for the overall success of the article,
  represented by demands that are $0.7, 1.0, 1.3$ times as large as a
  set of nominal demand values.
\end{itemize}
\markattention{This led to instances of ISPO with more than
  3\,500\,000 variables and constraints.}  In the following, we
present results on five such instances (results on all the other
instances we tried were almost the same):
\begin{itemize}
\item We measured for branch-and-bound the total CPU time in hours in
  the column denoted by ``\textsf{$t$[h}]''.

  Moreover, we counted how many
  \begin{itemize}
  \item exact computations of an ISPO with some prices fixed (see
    the column denoted by ``\textsf{\#ISPO (\%)}'')
  \item exact computations of LP relaxations of an ISPO with some
    prices fixed (see the column denoted by
    ``\textsf{\#ISPO$^\text{LP}$ (\%)}'')
  \end{itemize}
  we needed to find and prove an optimal solution.  In all other
  branch-and-bound nodes, it was sufficient to use the combinatorial
  bound coming from replacing the lot-type design restrictions to
  item-by-item supply.  The numbers in parentheses show the
  percentages of the numbers of all possible branch-and-bound nodes in
  order to indicate how often we got away with cheap bounds only.

  Moreover, we counted the number of exact computations of an ISPO
  until an optimal solution was found (but not yet proved) -- see the
  column denoted by ``\textsf{\#ISPO$^*$}''.  The column denoted by
  ``\textsf{$t^*$[h]}'' shows the CPU time in hours until this
  solution was found.
\item We measured for ping-pong the CPU time in minutes until no
  improvement happened anymore (see the column denoted by
  `\textsf{`$t$[min]}''). Moreover, we counted the number of
  iterations with improvements in the column denoted by
  ``\textsf{\#iter}''. Here, one iteration means one SOP and one POP
  computation. Finally, the column denoted by ``\textsf{Gap[\%]}''
  shows the relative optimality gap of the solution produced by
  ping-pong.
\end{itemize}

\begin{table}[ht!]
  \begin{center}\sffamily\footnotesize
    \begin{tabular}{r|*{5}{r}|*{3}{r}}
      \toprule
      \textbf{Instance} &
      \multicolumn{5}{c|}{\textbf{Branch\,\&\,Bound}} &
      \multicolumn{3}{c} {\textbf{Ping-Pong}}\\
      &$t$[h]&\#ISPO (\%)&\#ISPO$^\text{LP}$ (\%) &\#ISPO$^*$&$t^*$[h] &$t$[min]&\#iter & Gap[\%]\\
      \midrule
      P430204\_13&13.17&14 ($<$10$^{-\text{6}}$)&48 (1.17)&6&7.71&10.42&2&0.028\\  
      P430206\_13&15.31&24 ($<$10$^{-\text{4}}$)&31 (0.76)&1&0.32&3.77&1&0.015\\  
      P430207\_13&43.03&80 ($<$10$^{-\text{5}}$)&75 (1.83)&29&3.07&9.58&1.5&0.013\\   
      P490201\_13&60.78&45 ($<$10$^{-\text{5}}$)&177 (0.72)&23&31.67&12.62&2&0.023\\  
      P500206\_13&26.79&2 ($<$10$^{-\text{7}}$)&13 (1.19)&1&0.94&5.47&1.5&0.000\\  
      \midrule
      \textbf{$\varnothing$}&
      {\llap{\textbf{31.82}}}&
      {\llap{\textbf{33}}} \textbf{($<$10$^{-\text{5}}$)}&
      {\llap{\textbf{68.8}}} {\textbf{($1.13$)}}&
      {\llap{\textbf{12}}}&
      {\llap{\textbf{8.74}}}&
      {\llap{\textbf{8.37}}}&
      {\llap{\textbf{1.6}}}&
      {\llap{\textbf{0.016}}}\\
      \bottomrule
    \end{tabular}
  \end{center}

  \caption{Performance of the exact algorithm and the ping-pong heuristic}
  \label{tab:runtime-results}
\end{table}

The results provide evidence that
\begin{itemize}
\item the branch-and-bound algorithm can find and prove optimal
  solutions for production problems in a time that makes it suitable
  for benchmarking purposes; it is not fast enough for daily
  operation, because -- \markattention{even without any effort to prove
    optimality, the optimal solution is found too late};
\item the combinatorial dual-bound techniques help to avoid time
  consuming LP computations in many nodes;
\item the quality of ping-pong solutions is excellent;
\item the CPU times of ping-pong are in line with the real-time
  requirements of daily operation.
\end{itemize}
Thus, ping-pong could be routinely used in a field study designed as
in Section~\ref{sec_field_study}.

\subsection{Results of the field study}
\label{sec:fieldstudy}

Our test set of articles is denoted by $A$. For reasons of
comparability we consider for each branch the objective value of ISPO
divided per merchandise value over all articles from the set
$A$. We set the corresponding variables and parameters part
from for the different articles $a\in A$ by a
superscript~$a$. Apart from that the parameters name are identical to
the formulation of ISPO (Problem~\ref{sec:modelling}).

For each test-control pair of branches, the sums of objective function
values over all articles in $A$ were compared.  That means,
in particular, that expensive articles have a larger influence on the
result than cheap articles.  This point of view is in line with our
partner's point of view.

For reasons of comparability we consider the revenues measured by the
sums of objective values of ISPO divided by the maximal revenue
measured by the sums of merchandise values. This means for an initial
stock $\initialStock^a$ for the considered branch $b$, size $s$ and
article $a$ and a starting price $\pi_0^a$ we compute the
\emph{relative realized objective} of the independent non-anticipative
decisions

$\lotselect{} = \bigl(\lotselect{b,\ell,m}\bigr)_{b\in\Branches,
  \ell\in\Lottypes, m\in M}^a$ as
\begin{multline}
  \rro (x) = \frac{\text{objective achieved by $x$}}{\text{maximal
      possible objective}} = \\
  \frac{\displaystyle - \sum_{a\in\mathcal{A}} \sum_{\ell \in
      \Lottypes} \sum_{m \in \Multiplicities}
    \lotselect{b,\ell,m}^a \cdot \handlingcost_{b,\ell,m}^a -
    \sum_{i=1}^{\maxlots} \testaddlottypecost_i \cdot \lotcount{i}^a +
    \sum_{k \in \Periods} \exp(-\discount k) \Bigl( \sum_{a \in
      \mathcal{A}} \sum_{s \in \Sizes} \realizedyield{k,b,s}^a -
    \testmarkdowncost_k^a \realizedmarkdownused{k}^a \Bigr)}{- \sum_{a\in\mathcal{A}} \sum_{\ell \in
      \Lottypes} \sum_{m \in \Multiplicities}
    \lotselect{b,\ell,m}^a \cdot \handlingcost_{b,\ell,m}^a -
    \sum_{i=1}^{\maxlots} \testaddlottypecost_i \cdot \lotcount{i}^a + \sum_{a\in\mathcal{A}} \sum_{s \in \Sizes} \initialStock^a \cdot \pi_0^a}.
\end{multline}

Depending on article $a$, the entity $\lotcount{i}^a$ indicates that
an $i$th lot-type was used.  During the sales process, we observed
$\realizedyield{k,b,s}^a$ (the realized yield for branch $b$ and size
$s$ in period $k$) for article $a$ and $\realizedmarkdownused{k}^a$
(mark-down in period~$k$ -- yes or no).

Since we only consider a subset of branches we have to take into
account that pick costs, costs for additional lot types, and fixed
mark-down costs must be scaled with respect to the number of considered
branches. This way, we get a marginal cost $\testaddlottypecost_i$ for
the $i$th selected lot-type and mark-down costs
$\testmarkdowncost_k^a$ for period~$k$.  (For a complete notational
reference see Section~\ref{sec:modelling}).

The relative realized revenues are shown in Table~\ref{test_1} for
each pair at the second and third column. 

Wee see that on average, the use of the new
method gains almost two percentage points compared to the old method.

\begin{table}
 \centering \sffamily\footnotesize
\begin{tabular}{ccccc}
\toprule
\textbf{test-control-pair} & \textbf{$\rrotest$} & \textbf{$\rrocontrol$} & \textbf{$\rrotest-\rrocontrol$} & \textbf{signed rank}\\
\midrule
1&0.6333&0.6214&0.0119&2\\
2&0.6764&0.6080&0.0683&19\\
3&0.5919&0.6072&-0.0154&-5\\
4&0.6056&0.5898&0.0159&6\\
5&0.6637&0.5663&0.0974&26\\
6&0.6228&0.6031&0.0197&8\\
7&0.6377&0.6500&-0.0123&-3\\
8&0.5832&0.5845&-0.0013&-1\\
9&0.5968&0.5731&0.0237&11\\
10&0.5372&0.6276&-0.0904&-23\\
11&0.5651&0.5489&0.0163&7\\
12&0.5333&0.5904&-0.0571&-18\\
13&0.5782&0.5570&0.0212&9\\
14&0.6381&0.4940&0.1441&28\\
15&0.5054&0.5845&-0.0791&-21\\
16&0.5927&0.4993&0.0934&25\\
17&0.5872&0.4943&0.0929&24\\
18&0.6078&0.5691&0.0388&16\\
19&0.5762&0.6476&-0.0714&-20\\
20&0.5682&0.5323&0.0359&14\\
21&0.5133&0.4250&0.0883&22\\
22&0.5272&0.5547&-0.0275&-12\\
23&0.4015&0.5942&-0.1926&-30\\
24&0.4628&0.4860&-0.0232&-10\\
25&0.5168&0.4646&0.0522&17\\
26&0.5843&0.4621&0.1222&27\\
27&0.5658&0.4137&0.1521&29\\
28&0.4989&0.4608&0.0380&15\\
29&0.5466&0.5607&-0.0141&-4\\
30&0.5593&0.5272&0.0320&13\\
\midrule
\textbf{$\varnothing$} & 0.5692 & 0.5499 & 0.0193 &5.7\\
\bottomrule
\end{tabular}
\caption{RROs for the test-control-pairs -- all 81 articles}
\label{test_1}
\end{table}

In the following we use the controlled setup of the field study to
argue that the results are statistically significant with a prescribed
significance level of 5\,\% with no assumptions for the error
distributions (see, e.g., \cite{FreedmanPisaniPurves200702} for
general information on hypothesis testing).

We apply Wilcoxon signed-rank test from
statistics~\cite{Wilcoxon}. This test is applied for statistical
experiments for two related ordinal samples where no underlying
distribution can be assumed. It is an alternative to the Student's
t-test, which is applied for two related ordinal samples under the
assumption that the observations are normally distributed.

The procedure is as follows: The differences of the observations, here
$\rrotest-\rrocontrol$ -- at the fourth column of Table~\ref{test_1} are ordered increasingly according to their
absolute values. The ordering implies the corresponding rank for the
test-control-pair. Moreover, the sign of $\rrotest-\rrocontrol$ is
assigned to the rank. If the test branch won, than the rank has
positive sign, otherwise negative, see the fifth column of Table~\ref{test_1}. The rank sum is the sum of all
ranks with positive sign.  To check significance in terms of a better performance of the test branches, we compute the probability that
this or a higher rank-sum is observed by pure chance.

Our null-hypothesis is: using the new
method 
does not improve operations systematically.  That is, it does not
increase the probability to obtain a better objective function value
in practice.

The motivation for this test is: If the null-hypothesis is true, the
signed-rank sum would lead to a rank-sum close to~$\frac{n(n+1)}{4}$ with no
systematic positive deviation.\footnote{$\frac{n(n+1)}{2}$ is the sum of all ranks according to the Gaussian sum.} 

More specifically: With a predefined significance level of $\alpha$ we
can reject our null hypothesis ``test branches not systematically
better than control branches'' whenever we observe rank-sum $k$ and
$P_n(X\geq k) < \alpha$ where $n$ is the number of test-control-pairs.

For the data from Table~\ref{test_1} we get a rank-sum of 318. The
probability for getting an equal or higher rank-sum is $P_{30}(X\geq
318)\approx4.02\%$.
 
Thus, we can reject the null hypothesis with a significance level of
5\,\%.  Consequently, the test branches performed for the whole test
set of 81 articles significantly better than the control branches.

However, we could observe that there were some operational anomalities
like failed price cuts in the control branches. In order to estimate
the influence of the new method in the most conservative fashion, we
removed all articles which may have been affected by systematic
disturbances of operations.  This led to a second set of
articles~$A'$ with only 23 articles remaining.

The particular RROs are stated in Table~\ref{test_2}. We see that in
the case of heavily cleaned-up data the RRO for the test branches is
still more than $1.5$ percentage points higher than in the control
branches.  We repeated Wilcoxon signed-rank test for this smaller test
set.  Wilcoxon signed-rank test now yields a rank-sum of 271, which
leads to a probability of $P_{30}(X\geq 271)=22\%$ that a better
performance of the test-branches was observed by pure chance. Thus,
for the heavily cleaned-up data we still observe a relevant effect (1.5
percentage points improvement) whose observation can no longer be
testified as significant.  This is essentially caused by the fact that
for such a small (but relevant) effect the sample set $A'$
is simply no longer large enough to prove significance.  Still, the
probability for a systematic improvement is much larger than the
probability that the observed effect was caused by pure chance.

\begin{table}
 \centering
\sffamily\footnotesize
\begin{tabular}{ccccc}
\toprule
\textbf{test-control-pair} & \textbf{$\rrotest$} & \textbf{$\rrocontrol$} & \textbf{$\rrotest-\rrocontrol$} & \textbf{signed rank}\\
\midrule
1&0.4215&0.6673&-0.2458&-28\\
2&0.5874&0.4758&0.1116&18\\
3&0.6572&0.4865&0.1708&25\\
4&0.5948&0.4773&0.1175&21\\
5&0.5491&0.4153&0.1338&24\\
6&0.5799&0.5117&0.0682&13\\
7&0.4833&0.5454&-0.0621&-12\\
8&0.4648&0.5124&-0.0476&-9\\
9&0.5051&0.4923&0.0128&2\\
10&0.4933&0.6094&-0.1162&-19\\
11&0.4926&0.4998&-0.0071&-1\\
12&0.4205&0.4706&-0.0501&-10\\
13&0.4352&0.3746&0.0607&11\\
14&0.7046&0.2860&0.4186&30\\
15&0.4547&0.5281&-0.0734&-14\\
16&0.5146&0.3846&0.1300&22\\
17&0.5285&0.4247&0.1038&17\\
18&0.4802&0.5081&-0.0279&-3\\
19&0.3562&0.4865&-0.1303&-23\\
20&0.4119&0.4496&-0.0377&-5\\
21&0.2195&0.2577&-0.0382&-6\\
22&0.4274&0.5437&-0.1163&-20\\
23&0.2262&0.6415&-0.4153&-29\\
24&0.4006&0.3252&0.0754&16\\
25&0.3779&0.4244&-0.0465&-8\\
26&0.4759&0.4008&0.0750&15\\
27&0.5926&0.3971&0.1955&26\\
28&0.4458&0.4116&0.0342&4\\
29&0.4540&0.4985&-0.0445&-7\\
30&0.5278&0.3050&0.2228&27\\
\midrule
\textbf{$\varnothing$}&0.4761&0.4604&0.0157&2.57\\
\bottomrule
\end{tabular}
\caption{RROs for the test-control-pairs -- heavily cleaned-up data, 23 articles}
\label{test_2}
\end{table}
So far, we assessed the quality of the decisions of the various
methods on the basis of our objective function that was carefully
engineered together with our partner.  Yet, it is interesting to see
that the new two-stage method outperforms the old method in some very
important criteria at the same time.  In Table
\ref{kumulierte_Kennwerte} we list average \rro, relative gross
yields, and relative sales for all test-control-pairs.  
For both revenue and gross yield we see improvements by the new
method. In contrast to this, the number of sales is only minimally
smaller for the new method.

Now, which decisions have been taken differently by the new method?
On the heaviliy cleaned-up data set of 23 articles, the new price
optimization suggested alltogether 14 mark-downs in the test branches,
while the manual strategy in the control branches led to 18
mark-downs on the same set.  This difference may be caused by the fact
that the new method tries to balance the increase in sales against the
decrease in the yield per piece more thoroughly.

\begin{table}
 \label{kumulierte_Kennwerte}
  \begin{center}
    \begin{tabular}{cccc}
      \toprule 
      \textbf{sample} & \textbf{relative realized objective} & \textbf{gross yield} & \textbf{sales} \\
      \midrule
      test & 0.4761 & 0.6829 & 0.7951 \\
      control & 0.4604 & 0.6744 & 0.8021 \\
      \bottomrule
    \end{tabular}
  \end{center}
 \caption{Alternative performance metrics, heavily cleaned-up data.}
\end{table}

Table~\ref{tab:lottypesTestControl} shows the differences in the
lot-type designs of the new and the old method for the 23 remaining
articles.\footnote{Since the lot-type design of the control branches
  had to be reconstructed from in this respect incomplete transaction
  data, the multiplicities for the control branches do not always add
  up to~30. The lot-types are reliable, though.} The most obvious
effect is that the number of different lot-types used is usually
smaller for the new method than for the old method.  Since the old
method tries to approximate a fractional demand as closely as possible
by a supply distribution on the basis of suitable lot-types, it will
usually use as many lot-types as possible, even if the improvements of
a new lot-type are small. The goal of the new method is not to meet
the demand as closely as possible but to earn as much money as
possible.  Obviously, an additional lot-type is not always justified
by higher predicted profits in ISPO.  Consequently, ISPO does not
suggest to use such a new lot-type.  In the table we clearly see that
lot-type $(1, \ldots, 1)$ is very often used.  This is the result of
the rule that each branch has to receive at least one piece in every
size -- a fact that reduces the potential for improvement and should
be taken into account when the effect (1.5 to 2 percentage points
improvement) of using the new method is assessed.

\begin{table}[ht!]
  \begin{center}\sffamily\scriptsize
    \begin{tabular}{rll}
      \toprule
      \textbf{no.}&
      \textbf{lots delivered to test branches by new method}&\textbf{lots delivered
      to control branches by old method}\\
      \midrule
      1&4(2,2,3,4,3,3),19(1,1,1,1,1,1),7(1,1,2,2,2,2) &13(1,1,1,1,1,1),7(1,1,1,2,2,1),5(1,1,2,2,3,2),3(2,3,3,4,4,3) \\
      2&4(2,2,3,4,3,3),19(1,1,1,1,1,1),7(1,1,2,2,2,2) &13(1,1,1,1,1,1),7(1,1,1,2,2,1),5(1,1,2,2,3,2),3(2,3,3,4,4,3) \\
      3&5(2,2,3,4,3,3),18 (1,1,1,1,1,1),7(1,1,2,2,2,2) &15(1,1,1,1,1,1),7(1,1,1,2,2,1),3(2,3,3,4,4,3),3(1,1,2,2,3,2)\\
      4&5(2,2,3,4,3,3),18(1,1,1,1,1,1),7(1,1,2,2,2,2) &12(1,1,1,1,1,1),8(1,1,1,2,2,1),3(1,1,2,2,3,2),3(2,3,3,4,4,3) \\
      5&5(2,2,3,4,3,3),18(1,1,1,1,1,1),7(1,1,2,2,2,2) &13(1,1,1,1,1,1),7(1,1,1,2,2,1),4(1,1,2,2,3,2),4(2,3,3,4,4,3) \\
      6&5(2,2,3,4,3,3),18(1,1,1,1,1,1),7(1,1,2,2,2,2) &14(1,1,1,1,1,1),7(1,1,1,2,2,1),6(1,1,2,2,3,2),3(2,3,3,4,4,3) \\
      7&5(2,2,3,4,3,3),18(1,1,1,1,1,1),7(1,1,2,2,2,2) &14(1,1,1,1,1,1),7(1,1,1,2,2,1),6(1,1,2,2,3,2),3(2,3,3,4,4,3)\\
      8&5(2,2,3,4,3,3),18(1,1,1,1,1,1),7(1,1,2,2,2,2) &12(1,1,1,1,1,1),7(1,1,1,2,2,1),7(1,1,2,2,3,2),3(2,3,3,4,4,3) \\ 
      \midrule
      10&13(1,1,1,2,2,2),17(1,1,1,1,1,1) &6(2,2,2,3,4,4),8(1,1,1,2,3,3),12(1,1,1,2,2,2),4(1,1,1,1,1,1)\\
      11&13(1,1,1,2,2,2),17(1,1,1,1,1,1) &6(2,2,2,3,4,4),8(1,1,1,2,3,3),12(1,1,1,2,2,2),4(1,1,1,1,1,1) \\
      12&13(1,1,1,2,2,2),17(1,1,1,1,1,1) &6(2,2,2,3,4,4),8(1,1,1,2,3,3),12(1,1,1,2,2,2),4(1,1,1,1,1,1) \\
      14&13(1,1,1,2,2,2),17(1,1,1,1,1,1) &9(1,1,2,2,2,2),9(1,1,1,1,1,1),3(1,1,2,2,1,1),6(1,1,1,1,2,2) \\
      16&10(1,1,1,2,2,2),7(1,1,2,2,2,2),13(1,1,1,1,1,1) &14(1,1,2,2,3,3),5(2,2,3,4,4,4),11(1,1,1,2,2,2) \\
      17&10(1,1,1,2,2,2),7(1,1,2,2,2,2),13(1,1,1,1,1,1) &14(1,1,2,2,3,3),5(2,2,3,4,4,4),11(1,1,1,2,2,2) \\
      18&10(1,1,1,2,2,2),7(1,1,2,2,2,2),13(1,1,1,1,1,1) &14(1,1,2,2,3,3),5(2,2,3,4,4,4),11(1,1,1,2,2,2) \\
      \midrule
      19&18(3,2,1),12(2,1,1) &10(4,2,1),19(3,2,1) \\
      20&8(1,3,2),22(1,2,1) &10(1,2,1),6(2,4,3),11(1,3,2),2(1,1,1) \\
      21&8(1,3,2),22(1,2,1) &22(1,2,1),6(1,1,1),2(1,3,1) \\
      22&7(2,4,3),11(1,2,1),4(2,3,2),8(1,3,2) &16(1,2,1),7(2,4,3),3(1,3,2),1(1,2,2) \\
      23&18(3,2,1),12(2,1,1) &1(2,1,1),9(4,2,1),18(3,2,1),1(1,1,1) \\
      \bottomrule
   \end{tabular}
 \end{center}
 \caption{Supply for the test and control branches in terms of lots.}\label{tab:lottypesTestControl}
\end{table}

\begin{table}[ht!]
  \begin{center}\sffamily\footnotesize
    \begin{tabular}{ccccccc}
      \toprule
      &\textbf{com.~group}&\textbf{no.}&\textbf{article}&\textbf{predicted}&\textbf{realized}&\textbf{gap} \\
      \midrule
      &wof&1&1597437&527.96&318.64&-0.3965 \\
      &wof&2&1597438&527.96&285.88&-0.4585 \\
      &wof&3&1597791&900.39&490.29&-0.4555 \\
      &wof&4&1598323&900.39&603.34&-0.3299 \\
      &wof&5&1598324&900.39&533.66&-0.4073 \\
      &wof&6&1599002&391.20&482.88&0.2343 \\
      &wof&7&1599007&391.20&415.96&0.0633 \\
      &wof&8&1599843&700.42&656.63&-0.0625 \\
      &wof&9&1599850&700.42&521.56&-0.2554 \\
      \midrule
      &woc&10&1593027&620.24&440.91&-0.2891 \\
      &woc&11&1593028&620.24&435.63&-0.2976 \\
      &woc&12&1593029&620.24&497.15&-0.1984 \\
      &woc&13&1593055&957.87&666.63&-0.3041 \\
      &woc&14&1593056&957.87&545.82&-0.4302 \\
      &woc&15&1593057&957.87&622.68&-0.3499 \\
      &woc&16&1593079&631.84&680.98&0.0778 \\
      &woc&17&1593080&631.84&664.36&0.0515 \\
      &woc&18&1593081&631.84&651.06&0.0304 \\
      \midrule
      &wu&19&1595383&393.90&292.29&-0.2580 \\
      &wu&20&1597776&624.91&414.50&-0.3367 \\
      &wu&21&1597803&262.52&297.26&0.1324 \\
      &wu&22&1598044&371.35&508.94&0.3705 \\
      &wu&23&1599151&421.82&364.31&-0.1363 \\
      \hline
      $\mathbf{\varnothing}$& &&&&&-0.1742 \\
      \textbf{sd}& &&&&& 0.2388 \\
      \bottomrule
    \end{tabular}
  \end{center}

  \caption{Comparison of objective function values -- predicted by
    ISPO versus realized.}
  \label{tab:predictedOFV}
\end{table}

\begin{table}[ht!]
  \begin{center}\sffamily\footnotesize
    \begin{tabular}{ccccccc}
      \toprule
      &\textbf{com.~group}&\textbf{no.}&\textbf{article}&\textbf{predicted}&\textbf{realized}&\textbf{gap} \\
      \midrule
      &wof&1&1597437&232.295&194&-0.1649 \\
      &wof&2&1597438&232.295&177&-0.2380 \\
      &wof&3&1597791&231.724&180&-0.2232 \\
      &wof&4&1598323&231.724&210&-0.0937 \\
      &wof&5&1598324&231.724&198&-0.1455 \\
      &wof&6&1599002&227.202&235&0.0343 \\
      &wof&7&1599007&227.202&214&-0.0581 \\
      &wof&8&1599843&225.305&253&0.1229 \\
      &wof&9&1599850&225.305&226&0.0031 \\
      \midrule
      &woc&10&1593027&204.52&203&-0.0074 \\
      &woc&11&1593028&204.52&206&0.0072 \\
      &woc&12&1593029&204.52&207&0.0121 \\
      &woc&13&1593055&199.606&221&0.1072 \\
      &woc&14&1593056&199.606&204&0.0220 \\
      &woc&15&1593057&199.606&218&0.0922 \\
      &woc&16&1593079&225.639&251&0.1124 \\
      &woc&17&1593080&225.639&237&0.0504 \\
      &woc&18&1593081&225.639&235&0.0415 \\
      &wu&19&1595383&138.836&125&-0.0997 \\
      \midrule
      &wu&20&1597776&122.689&95&-0.2257 \\
      &wu&21&1597803&124.373&123&-0.0110 \\
      &wu&22&1598044&191.482&213&0.1124 \\
      &wu&23&1599151&138.836&130&-0.0636 \\
      \hline
      $\mathbf{\varnothing}$ & &&&&&-0.0267 \\
      \textbf{sd} & &&&&&0.1133 \\
      \bottomrule
    \end{tabular}
  \end{center}

  \caption{Comparison of sales -- predicted by ISPO versus realized.}
  \label{tab:predictedSALES}
\end{table}

In Tables~\ref{tab:predictedOFV} and \ref{tab:predictedSALES} we show
how well ISPO predicted the expected function values and the expected
sales, resp.  While the prediction quality of the expected function
values seems unsatisfactory, we get that the prediction of sales is
quite good.  That sales can be predicted well is more an indication
for the fact that essentially everything is sold anyway.  What matters
more is how much money can be earned by these sales.  And this in turn
indicates that it is vital to estimate the return when it comes to
deciding about the distribution of supply.  Although our predictions
are presumably biased (we usually predict better function values than
the realized ones), the volatility even in one commodity group is very
high (expressed by the standard deviation): a gap of zero is still
inside the interval ``average minus standard deviation'' through
``average plus standard deviation''. Yet, we will try to reduce the
bias of the prediction in the future by comparing realizations and
predictions more carefully.

\section{Conclusion and future work}
\label{sec_conclusion}

We introduced the Integrated Size and Price Optimization Problem ISPO,
which is a two-stage stochastic optimization problem with recourse to
optimize the distribution of goods among branches and sizes for a
fashion discounter. We presented an MILP formulation of the
deterministic equivalent in extensive form.  This model, however,
could not be solved for real-world instances by commercial MILP
software of the shelf. We therefore suggested one exact
branch-and-bound algorithm for benchmarking and a ping-pong heuristic
for daily production use.  In computational experiments on real-world
data we showed that the optimality gap of ping-pong is usually
tiny. In a five-month field study we applied ISPO in practice to
distribute produces over branches and sizes and observed the sales
process thereafter. We obtained an improvement for the realized
relative objective of more than $1.5$ percentage points compared to a
one-stage lot optimization model.  Because the field study was
designed as a controlled statistical experiment, we could show that
(for the complete set of the test articles) it is very unlikely that
an improvement happened by pure chance.

In order to be able to cope with more applicable lot-types, it would
be very interesting to generalize the branch-and-price algorithm for
the SLDP (size optimization only by solving the stochastic lot-type
design problem) in \cite{Kiessling+Kurz+Rambau:LDP-CG-Preprint:2012}
to the ISPO.  The ping-pong heuristic computes solutions to first and
second stage seperately with the variables of the other stage fixed;
thus, at least the ping-pong heuristic should also work with many
applicable lot-types.

Similarly, aggregating price selections to complete markdown
strategies like in Section~\ref{sec_algorithms} could possible be used
to generate a tighter ILP formulation with an independent decision
variable for each markdown strategy.  The more relevant advantage of
such a formulation, however, is the following.
Since our problem can be formulated as a two stage stochastic integer
linear program, one might apply corresponding general algorithms from
that area, see, e.g., \cite{1223.90001} for an overview. A quite
common algorithm is the so-called L-shaped method (stochastic Benders
decomposition of the second stage). A necessary condition for the
L-shaped method is, that the target function $Q(x,e)$ of the second
stage is concave and continuous (for maximization), which is often
violated in the presence of second-stage integrality constraints.
However, if variables for complete price trajectories are used instead
of periodical price assignment variables, then the second stage is
mathematically more well-behaved.
This might be a promising direction for further research.

The most important question posed by this work is, however: Can the
demand forecasts be improved by statistical methods, with which many
parameters can be estimated well by few observations?  This points
into the direction of support vector machines
\cite{Christmann+Steinwart:SVM:2008}. It would be interesting to
learn whether or not the practical impact of our optimization results
will improve if such more sophisticated forecasting methods are used
in practice.


\end{document}